\documentclass{amsart}

\usepackage{amssymb}
\usepackage[all]{xy}
\usepackage{hyperref}

\usepackage{enumitem}   

\setlist[enumerate]{itemsep=.2em,topsep=.2em,leftmargin=1.25em,itemindent=2.0em}

\newtheorem{thm}{Theorem}

\newtheorem{lem}[thm]{Lemma}

\theoremstyle{definition}

\newtheorem{say}[thm]{}
\newtheorem{exmp}[thm]{Example}

\newtheorem{rems}[thm]{Remarks} 

\newtheorem*{ack}{Acknowledgments}      
\newtheorem{notation}[thm]{Notation}   
  
\newtheorem{defn-thm}[thm]{Definition--Theorem}  
\newtheorem{defn-lem}[thm]{Definition--Lemma}  

\newtheorem{main-exmp}[thm]{Main Example}
\newtheorem{baby-exmp}[thm]{Baby Example}
\newtheorem*{obs-nn}{Observation}
\newtheorem*{assertion-nn}{Assertion}   

\theoremstyle{remark}

\renewcommand{\c}[0]{{\mathbb C}}  

\renewcommand{\o}[0]{{\mathcal O}} 

\renewcommand{\r}[0]{{\mathbb R}}

\newcommand{\sing}[0]{\operatorname{Sing}}

\def\loccoh#1.#2.#3.#4.{H^{#1}_{#2}(#3,#4)}

\DeclareMathAlphabet{\mathchanc}{OT1}{pzc}%
                                {m}{it}

\usepackage[all]{xy}\xyoption{dvips}

\begin{document}
\bibliographystyle{amsalpha}


 \title[Cartan's Theorem B for Real Analytic Sets] {A converse of Cartan's Theorem B\\ for Real Analytic Sets}
 \author{J\'anos Koll\'ar}

 \begin{abstract}
   In this lecture we prove a converse to  Cartan's Theorem~B for real analytic sets, due to Fernando and Ghiloni [Arxiv/2506.18347].
 \end{abstract}

  \maketitle

 Let $\Omega_{\r}\subset \r^n$ be an open subset, and
 $X_{\r}:= (f_1=\cdots=f_m=0) \subset\Omega_{\r}$ a closed, real analytic subset defined by the real analytic functions  $f_i$. These are called
 {\it C-analytic sets} by Whitney and Bruhat \cite{whi-bru}.

 A function $f$ on $X_{\r}$ is {\it real analytic} if  each $p\in X_{\r}$ has a neighborhood  $p\in B_p\subset \r^n$ such that $f|_{X_{\r}\cap B_p}$ 
is  also the restriction of a  real analytic function on $B_p$.

   $X_\r$ is {\it coherent} if each $p\in X_{\r}$ has a neighborhood
   $p\in B_p\subset \r^n$, and  real analytic functions $\{f_{p,i}\}$  on $B_p$ that generate the ideal sheaf of $X_\r\cap B_p$.

As a consequence of Cartan's Theorem~B,
if  $X_{\r}$ is coherent, then
 real analytic functions on $X$ extend to  real analytic functions on
 $\Omega_{\r}$; see  \cite[Sec.7(2)]{cartan-real}.
 
 Recently Fernando and Ghiloni proved the converse  \cite{fer-ghi}.
That is, if all real analytic functions on $X_{\r}$ extend to  
 $\Omega_{\r}$, then $X_{\r}$ is coherent and defined by  real analytic functions.
We give  a shortened version of their  proof.

As in \cite{cartan-real} and  \cite{whi-bru}, the plan is to reduce real analytic questions to the theory of complex analytic Stein spaces.

\begin{notation}\label{complex.say}
  The $f_i$ extend to some open neighborhood
  $\Omega_\r\subset \Omega\subset \c^n$, and define 
  a closed, complex analytic subset 
  $X\subset \Omega$. By \cite[Prop.1]{cartan-real} and \cite[Chap.3]{MR236418},  we may choose
  $\Omega$ to be 
    Stein, and such that
   $\Omega_\r= \Omega\cap \r^n$ and $X\cap \Omega_\r=X_\r$.
   We may also assume that $X$ is reduced, and $X_\r$ is Zariski dense in $X$;
   then $X$ is invariant under complex conjugation.
   We  think of $X$ as a germ around $X_\r$, and shrink  $\Omega$ whenever necessary.
   
   Let $\pi:Y\to X$ be the normalization, see
   \cite[Chap.VI]{MR217337}. Then $Y$ is also Stein
   \cite[Sec.V.1.1]{MR580152}, and it is  a germ around  $\pi^{-1}(X_\r)$.

 Let
 $Y^{\rm sm}_{\r}\subset Y$  be the set of smooth real points.
 It is a  Zariski dense subset of $Y$.

 Note that $Y^{\rm sm}_{\r}\subset \pi^{-1}(X_\r)$. The key is to understand when
 $Y^{\rm sm}_{\r}$ is Euclidean dense in $\pi^{-1}(X_\r)$.
 \end{notation}

 \begin{thm} \label{cartan.b.thm.2}
Using Notation~\ref{complex.say}, the following are equivalent
   \begin{enumerate}
     \item $Y^{\rm sm}_{\r}$ is Euclidean dense in $\pi^{-1}(X_{\r})$.
     \item  $X_\r$ is  coherent.
  \item Real analytic functions  on $X_\r$ extend to real analytic functions  on $\Omega_\r$.
   \end{enumerate}
 \end{thm}

 \begin{rems}
   Example~\ref{exmp.1} shows that a non-coherent set can be a topological manifold; see   \cite[Sec.11]{whi-bru}  for more examples, and \cite[Sec.3.E]{fer-ghi} for  pictures.

 We prove in Paragraphs~\ref{say.1}--\ref{lem.5}  that either of (\ref{cartan.b.thm.2}.2--3) implies  (\ref{cartan.b.thm.2}.1); this is taken from \cite{fer-ghi}.

 In Paragraph~\ref{1.to.2.3} we  show   that (\ref{cartan.b.thm.2}.1) implies  (\ref{cartan.b.thm.2}.2--3).  Note that the equivalence of (\ref{cartan.b.thm.2}.1) and (\ref{cartan.b.thm.2}.2) is essentially 
 \cite[Prop.12]{cartan-real}; see also \cite[p.94]{MR217337}.

 For the functions in Paragraph~\ref{lem.5} the  extension fails  at a single point. The arguments in \cite{fer-ghi}  produce many more examples, yielding an almost complete description of the possible non-extendability sets.

 \cite[1.6]{fer-ghi} proves the analog of Theorem~\ref{cartan.b.thm.2} for Nash subsets and Nash functions.
\end{rems}

 \begin{say}[Finding global equations] The following argument is in \cite[p.2]{fer-ghi}, which credits \cite[Thm.1]{MR1413268}.
   
As a preliminary construction, let $P\subset \Omega_\r$ be a discrete point set, disjoint from $X_\r$. If $f_1, \dots, f_r$ generate the ideal sheaf of $P$ in $\Omega$, then  $g_P:=\sum |f_i|^2$ 
   is a real analytic function on $\Omega_\r$ that  vanishes only on $P$. Let   $h_P$ be a real analytic function on $\Omega_\r$ that extends  $g^{-1}_P|_{X_{\r}}$. Then $F_P:=1-g_Ph_P$ is real analytic, vanishes on $X_\r$, but does not vanish on  $P$.  The $F_P$ are equations of $X_\r$.

   Now start with $F_{P_0}\equiv 0$. 
   Assume now that we already have equations $F_{P_0},\dots, F_{P_r}$.
   Let $P_{r+1}$ be any discrete point set that contains at least 1 point in each connected component of the smooth locus of  $(F_{P_0}=\cdots= F_{P_r}=0)\setminus X_\r$.
   After at most $m=1+\dim \Omega_\r$ steps, we get  that
   $X_\r=(F_{P_0}=\cdots= F_{P_m}=0)$. \qed
   \end{say}

 \begin{say}[Finding  non-coherence points]\label{say.1}
   Let $Y_i\subset Y$ the an irreducible component.
   Note that  $Y_i\cap \pi^{-1}(X_{\r})$ is connected, since
$Y$ is a germ around  $\pi^{-1}(X_\r)$.
   Thus if $Y^{\rm sm}_{\r}\subset \pi^{-1}(X_{\r})$ is not Euclidean dense, then  there is a $y_1\in Y$ that is a boundary point of $Y^{\rm sm}_{\r}\subset \pi^{-1}(X_{\r})$.

 Set $x_1:=\pi(y_1)$, choose $0<\epsilon\ll 1$,   let  $B(x_1,\epsilon)\subset \c^n$ denote the open ball of  radius $\epsilon$ with center $x_1$, and let 
 $W_j\subset \pi^{-1}\bigl(X\cap B(x_1,\epsilon)\bigr)$ be the
 irreducible components.  We may assume that, for every $j$,  
 \begin{enumerate}
 \item  $W_j$ contains a unique point  $y_j\in \pi^{-1}(x_1)$, 
 \item  $W_j\cap \pi^{-1}(X_{\r})$ is connected, and
   \item $y_1\in W_1$  is a boundary point of  $W_1\cap Y^{\rm sm}_{\r}\subset W_1\cap \pi^{-1}(X_{\r})$.
 \end{enumerate}
 Therefore  there is a  $y'_1\in W_1\cap\pi^{-1}(X_{\r})$ that is not in the closure of $Y^{\rm sm}_{\r}$.

 Set $x'_1=\pi(y'_1)$
 and choose  $0<\eta\ll \epsilon$ such that 
 the irreducible components  $V_j$ of 
$\pi^{-1}\bigl(X\cap B(x'_1,\eta)\bigr)$
 are in bijection with the points of $\pi^{-1}(x'_1)$.
 We index such that $y'_1\in V_1$.
 There are 2 possibilities:
 \begin{enumerate}\setcounter{enumi}{3}
   \item  either  $X_{\r}\cap B(x'_1,\eta)\subset \cup_{j\neq 1} \pi\bigl(V_j\cap \pi^{-1}(X_{\r})\bigr)$, 
\item  or we can perturb $y'_1$ a little to achieve that the 
  $V_j\cap \pi^{-1}(X_{\r})$ are empty for  $j\neq 1$.
   In this case $X_{\r}\cap B(x'_1,\eta)\subset \sing X$.
  (Thus $x'_1$ is on  a `tail' of  $X_{\r}$.)
 \end{enumerate}
 \end{say}

 \begin{say}[Proof of (\ref{cartan.b.thm.2}.2) $\Rightarrow$   (\ref{cartan.b.thm.2}.1)]

   By (\ref{say.1}.3)  real points are Zariski dense in $W_1$.
    Thus if  $g$ is analytic on
   $B(x_1,\epsilon)$, and vanishes along $X_\r\cap B(x_1,\epsilon)$,
   then $g\circ \pi$ vanishes along   $W_1$.

   However, the real points are not Zariski dense in
   $V_1\subset W_1$. Choose   $g'$  to be $1$ on $V_1$ and $0$ on the other $V_i$. With $H$ as in  Lemma~\ref{lem.6},
    $Hg'$ descends to an nonzero analytic function on
   $X\cap B(x'_1,\eta)$ that vanishes along $X_\r\cap B(x'_1,\eta)$.
   Thus coherence fails at $x_1\in X_\r$. \qed
 \end{say}

 In order to see  (\ref{cartan.b.thm.2}.3) $\Rightarrow$   (\ref{cartan.b.thm.2}.1),  we construct an analytic function $G$ on  $Y$
 such that $y'_1\in (G=0)\subset V_1$.
 Then we show that $HG^{-1}$ is real analytic on $X_\r$, but does not extend to $\Omega_\r$ for any such $G$.

 \begin{say}[Construction of $G$]\label{say.2}
    Set $x'_1=(p_1,\dots, p_n)$ and 
   start with the function $\phi:=\sum (x_i-p_i)^2$ and its zero divisor
   $D:=(\phi=0)\subset \Omega$. Note that $x'_1$ is the only real point of $D$, we may thus assume that
   $D\subset  B(x'_1,\eta)$  (possibly after shrinking $\Omega$).

    The pull-back  of $D$ to $Y$ decomposes as a disjoint union of  divisors
   $D_i\subset V_i$.
    We claim that   $\o_Y(-D_1)$ is a trivial line bundle. If this holds, then we can take $G$ to be any generating section of $\o_Y(-D_1)$.

    By the Oka principle, on a Stein space, a topologically trivial line bundle is also
    analytically trivial (cf.\  \cite[Sec.V.2.5]{MR580152}). Thus it is enough to construct a
     continuous generating section of $\o_Y(-D_1)$.

   Since $\phi$ is everywhere positive on the boundary of $\r^n\cap B(x'_1,\eta)$  there is a nowhere zero continuous function  $s$ on $X\cap B(x'_1,\eta)$ 
    that agrees with  $\phi$ on the  boundary  (again possibly after shrinking $\Omega$).

   We can now use the pull-back of $\phi$ outside  $\cup_{j\neq 1}V_j$,
   and the  pull-back of $s$ on  $\cup_{j\neq 1}V_j$
   to get the required   continuous generating section of $\o_Y(-D_1)$.
\end{say}

 \begin{say}[Step 1 of the proof of (\ref{cartan.b.thm.2}.3) $\Rightarrow$   (\ref{cartan.b.thm.2}.1)]\label{lem.4}
     We show that 
 $HG^{-1}$  descends to a real analytic function  $F$ on $X_\r$.

 If $x\in X$ is a smooth point, then set 
 $F(x):=HG^{-1}(\pi^{-1}(x))$.
 If $x\in X$ is a singular point, then we set $F(x)=0$. 

We need to check that  $F$ is real analytic.
Outside $B(x'_1,\eta)$ this follows from Lemma~\ref{lem.6}, since $G^{-1}$ is analytic on    $Y\setminus \pi^{-1}B(x'_1,\eta)$.

 Inside $B(x'_1,\eta)$, there are 2 possibilities.
 If we are in case (\ref{say.1}.5) then
 $X_{\r}\cap B(x'_1,\eta)\subset \sing X$, so $HG^{-1}$ is the 0 function on
 $X_{\r}\cap B(x'_1,\eta)$.

 In case (\ref{say.1}.4) consider
 $\cup_{j\neq 1}V_j\to X\cap B(x'_1,\eta)$.
 Its image  $X^{(1)}\subset X\cap B(x'_1,\eta)$ has one fewer irreducible components than  $X\cap B(x'_1,\eta)$, but $X^{(1)}_{\r}= X_{\r}\cap B(x'_1,\eta)$.

 Since $G^{-1}$ is analytic on $\cup_{j\neq 1}V_j$, 
  we can apply Lemma~\ref{lem.6} to
 $\cup_{j\neq 1}V_j\to X^{(1)}$, to get 
  an analytic function   $F^{(1)}$ on $X^{(1)}$,
which agrees with  $F$ on smooth points of $X^{(1)}$.
  The restriction of $F^{(1)}$ to
  $X^{(1)}_{\r}$ is real analytic. 
  \qed
\end{say}

 \begin{say}[Step 2 of the proof of (\ref{cartan.b.thm.2}.3) $\Rightarrow$   (\ref{cartan.b.thm.2}.1)]\label{lem.5} We show that 
 $F=HG^{-1}$   does not extend to a real analytic function on
   $\r^n\cap B(x_1,\epsilon)$.

  Assuming the contrary, there is an analytic extension $F_U$
 to a neighborhood $U$ of $\r^n\cap B(x_1,\epsilon)\subset \c^n$.
 Let  $F_1$  denote its pull-back to $W_1\cap \pi^{-1}(U)$.
 Note that  $F_1$ agrees with $HG^{-1}$ on
 real points.

 Let $ W_1^\circ\subset W_1\cap \pi^{-1}(U)$ be  the connected component that contains $y_1$.
 Since $y_1\in W_1$ is a boundary point of $Y^{\rm sm}_{\r}$, smooth real points of
 $ W_1^\circ$ are Zariski dense.  Thus
 $F_1$ agrees with $HG^{-1}$ on $ W_1^\circ$.
  By
 (\ref{say.1}.2) 
 $W_1\cap \pi^{-1}(X_{\r})$ is connected, so
 $y'_1\in W_1^\circ$.  Thus $F_1$ agrees with $HG^{-1}$ near $y'_1$. This is a contradiction since $HG^{-1}$ is not analytic at $y'_1$. \qed
\end{say}

 \begin{say}[Cartan's proof of (\ref{cartan.b.thm.2}.1) $\Rightarrow$   (\ref{cartan.b.thm.2}.2--3)]\label{1.to.2.3}
   Pick $x\in X_\r$ and let
   $Y_i\subset Y\cap \pi^{-1}B(x, \epsilon)$ be the irreducible components for some $0<\epsilon\ll 1$. By (\ref{cartan.b.thm.2}.1) and Lemma~\ref{complex.say.1}, the real points of $Y\cap \pi^{-1}B(x, \epsilon)$ are Zariski dense. Thus if $g$ is analytic on  $B(x, \epsilon)$ and vanishes on $X$, then $g\circ \pi$ also vanishes on $Y\cap \pi^{-1}B(x, \epsilon)$, hence  $g$  vanishes on $X\cap B(x, \epsilon)$.
   That is, the ideal sheaf of $X_{\r}$ is the same as the  
   ideal sheaf of the complex analytic $X$, which is coherent by Oka's theorem; see \cite[Chap.IV]{MR217337} or \cite[Sec.4.2]{MR580152}. 
   Thus (\ref{cartan.b.thm.2}.2) holds.

   Let next $g$ be a real analytic function on $X_\r$. By definition,
   $X$ is covered by complex balls  $B(x_i, \epsilon_i)$ such that
   the restriction on $g$ to $X_\r\cap  B(x_i, \epsilon_i)$ extends to a
   complex analytic $g^*_i$ on $B(x_i, \epsilon_i)$.
   Let $\bar g_i$ denote the restriction of $g^*_i$ to
   $X\cap B(x_i, \epsilon_i)$. We claim that the $\bar g_i$ define a
   complex analytic function $\bar g$ on $X$.

   If an overlap  $B(x_i, \epsilon_i)\cap B(x_j, \epsilon_j)$ is nonempty, then real points are Zariski dense in its preimage by $\pi$,
   thus the restrictions of $\bar g_i$ and of $\bar g_j$ agree.

   Finally $\bar g$ extends to a complex analytic function $g_{\Omega}$ on $\Omega$ by Cartan's Theorem~B on Stein spaces \cite{MR64154}.\qed

   \end{say}

 We used 2 straightforward lemmas.

 \begin{lem}\label{complex.say.1}  Let $U\subset \c^n$ be a complex space that is invariant under complex conjugation, and
   $p\in U$ a smooth real point. Then there is a unique irreducible component $U_p\subset U$ that contains $p$.  This $U_p$ is invariant under complex conjugation, and
   the real points are Zariski dense in it. \qed
 \end{lem}
 
 \begin{lem}\label{lem.6}  Let $X$ be a reduced Stein space with normalization $\pi:Y\to X$.  Then  $\{h: h\cdot \pi_*\o_Y\subset \o_X\}$
   is a coherent ideal sheaf. By Cartan's Theorem~A there is  
   an  analytic function $H$ on $X$  such that
   $H\cdot \pi_*\o_Y\subset \o_X$, but $H$  does not vanish on  any irreducible component of $X$.
 Note that  $H$ vanishes on $\sing X$. \qed
   \end{lem}

 \begin{exmp} \label{exmp.1}  The traditional non-coherent example is the Whitney umbrella
   $(z^2=x^2y)$, but the following modified version is more subtle.
   
   Let $X_\r:=(z^3=x^3y)\subset \r^3$. Note that $z=x\sqrt[3]{y}$, hence $X_\r$ is  homeomorphic to $\r^2$. It is not coherent at the origin. The ideal at the origin is generated by $z^3-x^3y$.
   On the set $(y\neq 0)$  we also have $z-x\sqrt[3]{y}$, which is not a multiple of $z^3-x^3y$.
   \end{exmp}

\begin{ack}  I thank  J.F.~Fernando and R.~Ghiloni
  for  numerous comments and references.
   Partial  financial support    was provided  by the Simons Foundation   (grant number SFI-MPS-MOV-00006719-02).
    \end{ack}


\begin{thebibliography}{WB59}

\bibitem[Car53]{MR64154}
Henri Cartan, \emph{Vari\'et\'es analytiques complexes et cohomologie},
  Colloque sur les fonctions de plusieurs variables, tenu \`a{} {B}ruxelles,
  1953, Georges Thone, Li\`ege, 1953, pp.~41--55. \MR{64154}

\bibitem[Car57]{cartan-real}
\bysame, \emph{Vari\'et\'es analytiques r\'eelles et vari\'et\'es analytiques
  complexes}, Bull. Soc. Math. France \textbf{85} (1957), 77--99. \MR{94830}

\bibitem[FG25]{fer-ghi}
José~F. Fernando and Riccardo Ghiloni, \emph{A converse to {C}artan's
  {T}heorem {B}: The extension property for real analytic and {N}ash sets},
  \url{https://arxiv.org/abs/2506.18347}, 2025.

\bibitem[GR79]{MR580152}
Hans Grauert and Reinhold Remmert, \emph{Theory of {S}tein spaces}, Grundlehren
  der Mathematischen Wissenschaften, vol. 236, Springer-Verlag, Berlin-New
  York, 1979, Translated from the German by Alan Huckleberry. \MR{580152}

\bibitem[Nar66]{MR217337}
Raghavan Narasimhan, \emph{Introduction to the theory of analytic spaces},
  Lecture Notes in Mathematics, vol. No. 25, Springer-Verlag, Berlin-New York,
  1966. \MR{217337}

\bibitem[NT96]{MR1413268}
G.~Nardelli and A.~Tancredi, \emph{A note on the extension of analytic
  functions off real analytic subsets}, Rev. Mat. Univ. Complut. Madrid
  \textbf{9} (1996), no.~1, 85--97. \MR{1413268}

\bibitem[Tog67]{MR236418}
Alberto Tognoli, \emph{Propriet\`a{} globali degli spazi analitici reali}, Ann.
  Mat. Pura Appl. (4) \textbf{75} (1967), 143--218. \MR{236418}

\bibitem[WB59]{whi-bru}
H.~Whitney and F.~Bruhat, \emph{Quelques propri\'et\'es fondamentales des
  ensembles analytiques-r\'eels}, Comment. Math. Helv. \textbf{33} (1959),
  132--160. \MR{102094}

\end{thebibliography}

\def\cprime{$'$} \def\cprime{$'$} \def\cprime{$'$} \def\cprime{$'$}
  \def\cprime{$'$} \def\dbar{\leavevmode\hbox to 0pt{\hskip.2ex
  \accent"16\hss}d} \def\cprime{$'$} \def\cprime{$'$}
  \def\polhk#1{\setbox0=\hbox{#1}{\ooalign{\hidewidth
  \lower1.5ex\hbox{`}\hidewidth\crcr\unhbox0}}} \def\cprime{$'$}
  \def\cprime{$'$} \def\cprime{$'$} \def\cprime{$'$}
  \def\polhk#1{\setbox0=\hbox{#1}{\ooalign{\hidewidth
  \lower1.5ex\hbox{`}\hidewidth\crcr\unhbox0}}} \def\cdprime{$''$}
  \def\cprime{$'$} \def\cprime{$'$} \def\cprime{$'$} \def\cprime{$'$}
\providecommand{\bysame}{\leavevmode\hbox to3em{\hrulefill}\thinspace}
\providecommand{\MR}{\relax\ifhmode\unskip\space\fi MR }
\providecommand{\MRhref}[2]{%
  \href{http://www.ams.org/mathscinet-getitem?mr=#1}{#2}
}
\providecommand{\href}[2]{#2}

 \bigskip

  Princeton University, Princeton NJ 08544-1000, \

  \email{kollar@math.princeton.edu}

\end{document}